\documentclass[12pt]{article}

\usepackage{amscd}

\begin{document}

\title{Reciprocity law for compatible systems of abelian mod $p$ Galois representations}

\author{Chandrashekhar Khare}

\date{}

\maketitle

\newtheorem{theorem}{Theorem}
\newtheorem{lemma}{Lemma}
\newtheorem{prop}{Proposition}
\newtheorem{cor}{Corollary}
\newtheorem{example}{Example}
\newtheorem{conjecture}{Conjecture}
\newtheorem{definition}{Definition}
\newtheorem{quest}{Question}
\newtheorem{memo}{Memo}
\newcommand{\rhobar}{\overline{\rho}}
\newcommand{\Sha}{{\rm III}}
\newtheorem{conj}{Conjecture}

\noindent{\bf Abstract:} The main result of the paper 
is a {\bf reciprocity law} which proves that
compatible systems of semisimple, abelian mod $p$ representations (of arbitrary dimension)
of absolute Galois groups of number fields, arise from Hecke characters.
In the last section analogs for Galois groups of function fields of these
results are explored, and a question is raised whose answer will
require developments in transcendence theory in characteristic $p$.

\section{Introduction}

Motives defined over number fields 
give rise to compatible systems of $p$-adic and
mod $p$ representations of absolute Galois groups 
of number fields. Compatible systems of 
$p$-adic representations have been extensively studied.
Although compatible systems of mod
$p$ representations have received less attention, 
they have been considered by Serre in his
work in the 1960's and 1970's on studying adelic images of Galois groups acting
on products of $p$-adic Tate modules of an elliptic curve
for varying $p$.

In [K] we showed how studying compatible systems of mod $p$
representations can be useful in proving results about $p$-adic representations.
Specifically, using this point 
of view, we rederived in a simple way 
the result (see [He]) that strictly compatible systems 
of one-dimensional $p$-adic representations (see I-11 of [S]) arise from
Hecke characters. Compatible systems of mod $p$ representations 
make quite apparent how to use the fact that 
we have a compatible system rather than just one representation at hand,
while in the case of compatible $p$-adic systems this is not so apparent.
The difference between these two types of compatible systems
is mainly accounted for by the fact that given a $p$-adic representation 
it can be made part of at most one 
(semisimple) compatible system, 
while this is far from being true for a mod $p$ representation.

In [K] we proved that one-dimensional compatible mod $p$ systems
arise from Hecke characters. In this paper we would like
to generalise results of [K] to the case of higher dimensional, but still 
abelian, compatible mod $p$ systems. 
We recall the definition of compatible systems of $n$-dimensional, 
mod $p$ representations of the absolute Galois group of a number field $K$.
\begin{definition}\label{compatible}
 Let $K$ and $L$ be number fields and $S,T$ finite sets of places 
 of $K$ and $L$ respectively. An $L$-rational 
 strictly compatible system $\{\rho_{\wp}\}$ 
 of $n$-dimensional mod $\wp$ representations of
 $G_K:={\rm Gal}(\overline{K}/K)$ with defect set $T$ and ramification set $S$,
 consists of giving for each finite place $\wp$ of $L$ not in $T$ a 
 continuous, semisimple representation $$\rho_{\wp}:G_K \rightarrow GL_n({\bf
 F}_{\wp}),$$ for ${\bf F}_{\wp}$ the residue field of ${\cal O}_L$ at
 $\wp$ of characteristic $p$, that is \begin{itemize} \item  unramified at the places
 outside $S \cup$ $\{$ places of $K$
 above $p$ $\}$ \item  for each place 
 $r$ of $K$ not in $S$ there is a monic polynomial $f_r(X) \in L[X]$ such that for all places $\wp$ of
 $L$ not in $T$, coprime to the residue characteristic of $r$, and
 such that $f_r(X)$ has coefficients that are integral at $\wp$, 
 the characteristic polynomial of $\rho_{\wp}({\rm Frob}_r)$
 is the reduction of $f_r(X)$ mod $\wp$, where ${\rm Frob}_r$ is the
 conjugacy class of the Frobenius at $r$ in the Galois
 group of the extension of $K$ that is the fixed field of the kernel of
  $\rho_{\wp}$. \end{itemize} 
 \end{definition}
The following theorem was proved in [K]. 
\begin{theorem}
 An $L$-rational strictly compatible system $\{\rho_{\wp}\}$ 
 of one-dimensional mod $\wp$ representations of
 ${\rm Gal}(\overline{K}/K)$ arises from a Hecke character.
\end{theorem}
In Section 4 of [K] it is explained what one means by
saying that a compatible system of 1-dimensional mod $\wp$ representations
arises from a Hecke character. 

The theorem below is the main result of this paper. It generalises the result of [K] to
{\it semisimple} compatible systems $(\rho_{\wp})$: by this we mean
that $(\rho_{\wp})$ is a strictly compatible system as in the definition above 
and each $\rho_{\wp}$ has abelian image. By saying that a compatible system
$(\rho_{\wp})$ as in the definition above, is the sum of compatible systems $(\rho_{i,\wp})$ ($i=1,\cdots,n$),
we mean that $\rho_{\wp} \simeq \oplus_{i=1}^n\rho_{i,\wp}$ for all $\wp$ outside a finite defect set.
\begin{theorem}\label{main}
 An $L$-rational strictly compatible system $\{\rho_{\wp}\}$ 
 of abelian, semisimple mod $\wp$ representations of
 ${\rm Gal}(\overline{K}/K)$ is a direct sum of one-dimensional 
 compatible systems
 each of which arises from a Hecke character.
\end{theorem}
This theorem is but a small step in studying
compatible systems of mod $p$ Galois representations of arbitrary 
dimensions. Because of it such compatible systems which are abelian
are completely understood. 
Theorem \ref{main} also proves Conjectures 1 and 2 of [K] for all abelian compatible 
systems of mod $\wp$ Galois representations which might be summarised by saying that
abelian compatible mod $p$ systems of Galois representations are motivic 
(and all that this entails!).

The proof of Theorem \ref{main} uses the ideas of [K].
The main problem that arises when generalising the arguments in [K] is that although using arguments of [K]
one can easily get that the characteristic polynomials $f_r(X)$ each have a shape that is consistent with
the compatible system arising from Hecke characters, we cannot use the arguments in [K]
directly to show that the nature of the roots of $f_r(X)$ for varying $r$ is consistent with
the compatible system arising from Hecke characters.
It is this difficulty that is overcome in the lemma and corollary below and the ensuing arguments.
We focus mainly on this difficulty as otherwise the arguments are 
similar to [K] where the theorem above was announced.
 Theorem \ref{main} would have directly followed from Theorem 1 of [K] if we could prove
that the compatible system in the theorem is the sum of 1-dimensional compatible systems.
But this we  {\bf cannot} prove {\it a priori}. More generally it seems hard to prove that a compatible
system of mod $\wp$ or $\wp$-adic representations 
is ``the sum of 2 compatible systems'' if each mod $\wp$ or $\wp$-adic representation in
the compatible system is decomposable.

As in [K], Theorem \ref{main} has the following corollary.

\begin{cor}\label{taniyama}

 1. An $L$-rational strictly compatible system $\{\rho_{\wp}\}$ 
 of abelian, semisimple mod $\wp$ representations of
 ${\rm Gal}(\overline{K}/K)$ lifts to a compatible system of $n$-dimensional 
  $\wp$-adic representations
  $\{\rho_{\wp,\infty}\}$ as in I-11 of [S].

 2. An $L$-rational strictly compatible system of abelian, semisinmple
  $n$-dimensional 
  $\wp$-adic representations
  $\{\rho_{\wp,\infty}\}$ as in I-11 of [S]  
  arises from Hecke characters.
\end{cor}
The seocnd part of the corollary is implied by difficult results of Waldschmidt in transcendenatl number theory, while the proof we offeris more ``elementary''.

In the last section of the paper we indicate what the analogs of our theorems for function fields
should be. We do not prove these proposed analogs. We hope that someone will
carry out the proofs of these analogs and more interestingly answer the question 
at the end that we do know how to answer. The question asks for the natural analog
to a result of [He] in the setting of function fields: to answer it will probably need
transcendence results in characteristic $p$.

\section{Proof of theorem}

Let $(\rho_{\wp})$ be an abelian, semisimple, compatible system
of $L$-rational representations of $G_K$ with $K$ a number field,
that has finite defect set $T$ and 
finite ramification set $S$.
Using the arguments in Lemma 1 of [K], which generalise easily to higher dimensional abelain semisimple
compatible systems of mod $\wp$ representations, we can reduce to the case when
$K$ is a Galois extension of ${\bf Q}$ that contains $\sqrt{-1}$, $L$ contains
$K$ and is again Galois over ${\bf Q}$.
 
Using class field theory as in the proof of Theorem 1 of [K], the
association of Galois representations to Hecke characters
due to Taniyama and Weil that is recalled in Section 4.1 of [K],
the analog of Lemma 1 of [K] in this setting, and Proposition 2 in [K]
we are reduced to looking at a system of homomorphisms
$\rho_{\wp}:{\rm Cl}_{{\sf m}_{\wp}p} \rightarrow GL_n({\bf F}_{\wp})$.
Here ${\rm Cl}_{{\sf m}_{\wp}p}$ is the strict ray class group
of conductor ${\sf m}_{\wp}p$, with ${\sf m}_{\wp}$ the prime to $p$
part of the Artin conductor of $\rho_{\wp}$
which is divisible only by the primes in $S$.
By using the fact that the subgroup of ${\rm Cl}_{{\sf m}_{\wp}p}$ 
that is the image of principal ideals prime to ${\sf m}_{\wp}p$
is of index that is bounded independently of $\wp$, 
to considering the induced system of homomorphisms 
$({\cal O}_K/{\sf m}_{\wp}p{\cal O}_K)^* \rightarrow GL_n({\bf F}_{\wp})$
which by abuse of notation we denote by the same symbol $\rho_{\wp}$:
we shall now exclusively study only these homomorphisms. 
We have to show how these homomorphisms {\it arise} from algebraic characters
of $K^*$, where $K^*$ is considered as the ${\bf Q}$-valued points 
of the algebraic group ${\bf Res}_{K/{\bf Q}}({\bf G}_m)$ over ${\bf Q}$, 
that have a subgroup of finite index of the units ${\cal O}_K^*$ in their 
kernel (see Section 4.1 of [K]). 

These homomorphisms all
factor through the quotient by the image of the global units $E$ of $K$.
The compatible data translates into saying that for any principal
prime ideal generated by  an element $r$ of $K$ (we denote the ideal $(r)$ generated by $r$, again by $r$) 
that is not in $S$ there is a monic polynomial
$f_r(X) \in L[X]$  such that for all $\wp$ not in $T$ an whose residue characteristic is
different from that of $r$ and at which $f_r(X)$ is integral, the characteristic polynomial of $\rho_{\wp}({\rm Frob}_{r})$
is the reduction mod $\wp$ of $f_r(X)$.
Let the roots of $f_r(X)$ be $\gamma_{1,r},\cdots, \gamma_{n,r}$: note that when
$r$ varies all these are algebraic integers which lie in extensions
of $L$ that are of degree over $L$ that is bounded by $n$. By using Proposition 3 of [K], which
generalises a result of [CS], and using the compatible data
it is easy to see (see Section 4.2 of [K]) that there is an integer $t_{i,r}$ such that
$\gamma_{i,r}^{t_i,r}=\Pi_{\sigma}\sigma(r)^{m'_{i,r,\sigma}}$, $i=1,\cdots,n$,
with the exponents some integers and with $\sigma$ running through the distinct
embeddings of $K$ in $\overline{\bf Q}$. Note that the ramification
indices of all primes in any field $M$ of degree over ${\bf Q}$ bounded by some integer $T$, 
is bounded independently of $M$, and the number of roots
of unity in such $M$ is also bounded independently of $M$. From this it is easy to see that
there is an integer $m$ independent of $i$ and $r$ such that $\gamma_{i,r}^m=\Pi_{\sigma}\sigma(r)^{m_{i,r,\sigma}}$, $i=1,\cdots,n$,
with the exponents some integers.

We would like to prove that the ${m_{i,r,\sigma}}$'s are independent of choice of 
the principal prime ideal $(r)$ that is coprime to the primes in $S$. To prove this we 
cannot use the argument in [K] that relied on the fact that in the one-dimensional situation
the characteristic polynomial of $\rho_{\wp}(rr')$ (with $rr'$ regarded as the image of $rr'$ mod ${\sf m}_{\wp}p$)
is the reduction mod $\wp$ of  a fixed polynomial in $L[X]$ for almost all primes $\wp$ and where $r$ and $r'$
generate principal prime ideals, not in $S$. In the higher dimensional case  this is not {\it a priori}
the case, although {\it a posteriori} we will see that this is true. Thus we need to modify the arguments 
of [K] to take into account this complication, and indeed this is the main contribution of the present paper.

Consider the homomorphisms
$\rho_{\wp}^m:({\cal O}_K/{\sf m}_{\wp}p{\cal O}_K)^* \rightarrow GL_n({\bf F}_{\wp})$ 
(where by $mth$ power we just mean taking $mth$ powers of the $n$ homomorphisms that consititute $\rho_{\wp}$).
Choose any principal prime ideals $(r)$ and $(r')$ of $K$ that are not in $S$. 
Fix $i$ between 1 and $n$. Then we see that for almost all primes  $\wp$ there is a $i(\wp)$
that lies between 1 and $n$ such that the reduction of $\Pi_{\sigma}\sigma(r)^{m_{i,r,\sigma}}
\Pi_{\sigma}\sigma(r')^{m_{i(\wp),r',\sigma}}$ mod $\wp$ is a root of the characteristic 
polynomial of $\rho_{\wp}^m(rr')$. 

We have the following lemma:

\begin{lemma}
Let $(r)$ and $(r')$ be principal prime ideals of $K$ not in $S$. Fix an integer $i$ between 1 and $n$. 
Consider a prime $\ell$ of ${\bf Q}$ that is prime to the residue characteristics
of the primes in $S$ and is prime to the cardinalities of the multiplicative groups
of the residue fieds at primes of $S$. Then if a prime of $K$ splits completely
in $K(\zeta_{\ell},(\sigma(rr'))^{1/\ell})$ with $\sigma$ running through ${\rm Gal}(K/{\bf Q})$, it also splits completely in
one of the fields  $K(\zeta_{\ell},(\Pi_{\sigma}\sigma(r)^{m_{i,r,\sigma}}
\Pi_{\sigma}\sigma(r')^{m_{j,r',\sigma}})^{1/\ell})$ for some $j$ between 1 and $n$.
\end{lemma}

\noindent{\bf Proof:} Fix $i$ between 1 and $n$. Then we see that for almost all primes $\wp$ there is a $i(\wp)$
that lies between 1 and $n$ such that the reduction of $\Pi_{\sigma}\sigma(r)^{m_{i,r,\sigma}}
\Pi_{\sigma}\sigma(r')^{m_{i(\wp),r',\sigma}}$ mod $\wp$ is a root of the characteristic 
polynomial of $\rho_{\wp}^m(rr')$. Now for almost all primes $s$ of $K$ that split completely
in $K(\zeta_{\ell},(\sigma(rr'))^{1/\ell})$ for all $\sigma \in {\rm Gal}(K/{\bf Q})$, the cardinality of the residue field at
$s$ is 1 mod $\ell$ and further $\sigma(rr')$ is a $\ell$th power modulo $s$ for all $\sigma \in {\rm Gal}(K/{\bf Q})$.
Thus by choice of $\ell$, the reduction of $\Pi_{\sigma}\sigma(r)^{m_{i,r,\sigma}}
\Pi_{\sigma}\sigma(r')^{m_{i(\wp),r',\sigma}}$ mod $s$ is an $\ell$th power, and thus $s$ also splits 
in $K(\zeta_{\ell},(\Pi_{\sigma}\sigma(r)^{m_{i,r,\sigma}}
\Pi_{\sigma}\sigma(r')^{m_{i(\wp),r',\sigma}})^{1/\ell})$ which proves the lemma.

\begin{cor}
 Fix an integer $i$ between 1 and $n$. 
 For all sufficiently large primes $\ell$ of ${\bf Q}$, the subgroup generated by
 $\tau(rr')$ of $K^*/(K^{*})^{\ell}$ where $\tau$ runs through ${\rm Gal}(K/{\bf Q})$,
 contains the image of $\Pi_{\sigma}\sigma(r)^{m_{i,r,\sigma}}
 \Pi_{\sigma}\sigma(r')^{m_{j,r',\sigma}}$ in $K^*/(K^{*})^{\ell}$ for some $j$ between 1 and $n$.
\end{cor}

\noindent{\bf Proof:} Consider $K^*/(K^{*})^{\ell}$ as a ${\bf F}_{\ell}$ vector-space, and
denote the  ${\bf F}_{\ell}$ vector-space generated by
 $\sigma(rr')$ of $K^*/(K^{*})^{\ell}$ for $\sigma \in {\rm Gal}(K/{\bf Q})$ by $V$.
Let $W_j$ be the one-dimensional vector space of $K^*/(K^{*})^{\ell}$ generated by $\Pi_{\sigma}\sigma(r)^{m_{i,r,\sigma}}
\Pi_{\sigma}\sigma(r')^{m_{j,r',\sigma}}$ for $j$ between 1 and $n$. Let ${\bf V}$ be the span
of the $V$'s and the $W_j$'s. Then using Kummer theory, the Cebotarev density theorem, the above
lemma and the injectivity of the map $K^*/(K^*)^{\ell} \rightarrow K(\zeta_{\ell})^*/(K(\zeta_{\ell})^*)^{\ell}$
(note that $\sqrt{-1} \in K$: for this injectivity  see Lemma 2.1 of [CS])
we conclude that if an element of the dual ${\bf V}^*$ has kernel that contains $V$, then its
kernel also contains $W_j$ for some $j$ between 1 and $n$. But if $\ell$ is large enough, this forces 
one of the $W_j$'s to be contained in $V$ by the following easy claim.

\noindent{\bf Claim:} If a prime $\ell$ is bigger than a given integer $k$, 
any finite dimensional vector space over ${\bf F}_{\ell}$ cannot be written as the
union of $k$ proper subspaces. 

We apply this claim to the vector space $X=({\bf V}/V)^*$ 
and the subspaces $(V/W_j)^* \cap X$ where the intersection is taking place in ${\bf V}^*$.
From the claim we see that $W_j$ is contained in $V$ for some $j$ and thus the corollary follows.

\vspace{3mm}

We now claim that that fixing a prime $(r)$ not in $S$ and which lies above a prime of ${\bf Q}$ that splits completely in $K$
(we call such a prime a split prime),
for all split primes
$(r')$ of $K$ of different residue characteristic from that of $r$, 
the distinct tuples that occur in $\langle (m_{i,r,\sigma})_{\sigma} \rangle$
and $\langle (m_{i,r',\sigma})_{\sigma} \rangle $ ($i=1,\cdots,n$)
are the same. Namely, fixing an $i$ between 1 and $n$, from the corollary it follows that
there is a fixed $j(i)$ between 1 and $n$, such that for infinitely many primes $\ell$ 
the image of $\Pi_{\sigma}r^{m_{i,r,\sigma}}\Pi_{\sigma}r'^{m_{j(i),r',\sigma}} $
in $K^*/(K^{*})^{\ell}$ is contained in the subgroup generated by the images
of $\sigma(rr')$ as $\sigma$ runs through ${\rm Gal}(K/{\bf Q})$. From this as in Step 3 of proof
of Theorem 1 of [CS], using the unit theorem, 
we see that some power of $\Pi_{\sigma}r^{m_{i,r,\sigma}}\Pi_{\sigma}r'^{m_{j(i),r',\sigma}} $
is contained in the subgroup of $K^*$ generated by $\sigma(rr')$ as $\sigma$ runs through ${\rm Gal}(K/{\bf Q})$.
From this, as $r$ and $r'$ generate {\it split} primes of different residue characteristic, we see
that $m_{i,r,\sigma}=m_{j(i),r',\sigma}$. Thus we get an injection from the 
distinct tuples in the collection $\langle (m_{i,r,\sigma})_{\sigma} \rangle$ to the distinct tuples in
the collection $\langle (m_{i,r',\sigma})_{\sigma} \rangle $. Now as $r,r'$ play symmetric roles the claim follows.
Repeating the argument by fixing a split principal prime of residue charateristic prime to $r$
we in fact conclude that the distinct tuples that occur in the collection
$\langle (m_{i,r,\sigma})_{\sigma}\rangle$ is independent of $r$ with $r$ any split principal
ideal prime to $S$. 

We now need a small 
argument to prove in fact that even 
the multiplicities with which distinct tuples occur in $\langle (m_{i,r,\sigma})_{\sigma} \rangle$
is independent of $r$. 
The difficulty here is related to the fact that if we have a linear representation
$\rho$ of a group $G$ such that for every $g \in G$, $\rho(g)$ has 1 as an eigenvalue,
it does not follow that the identity representation occurs in the semisimplification of $\rho$. 

We know by our work that there are only
finitely many possibilities for the $m_{i,r,\sigma}$'s as $i,r,\sigma$  vary
($r$ varies over split principal
prime ideals not in $S$)
and let $N$ be a natural number greater than all the finitely many integers $|2m_{i,r,\sigma}|$'s.
Let $(\alpha)$ be a split prime ideal of $K$ not in $S$
and choose a large enough rational prime $p'$, coprime to the places in $S$ and $T$ and $\alpha$, such that whenever
$\Pi_\sigma \sigma(\alpha)^{m_{\sigma}}-1$ is not coprime to $p'$, 
with $m_{\sigma}$ integers 
and $|m_{\sigma}| \leq N$, then all the $m_{\sigma}$'s are 0.
This is possible as the $\sigma(\alpha)$'s for $\sigma \in {\rm Gal}(K/{\bf Q})$ are multiplicatively independent.
Consider integral elements $\beta$ in $K$ such that $\beta$ is congruent to $\alpha$ mod $p'$,
and $\beta$ generates a split prime ideal not in $S$.
We claim that the unordered collection of $n$, $[K:{\bf Q}]$-tuples
$\langle (m_{i,\alpha,\sigma})_{\sigma} \rangle$ for
$i=1,\cdots,n$, is the same as 
$\langle (m_{i,\beta,\sigma})_{\sigma} \rangle$ for
$i=1,\cdots,n$. This is because
the numbers $\Pi_{\sigma \in {\rm Gal}(K/{\bf Q})}\sigma(\alpha)^{m_{i,\alpha,\sigma}}$
and $\Pi_{\sigma \in {\rm Gal}(K/{\bf Q})}\sigma(\beta)^{m_{i,\beta,\sigma}}$,
with the latter congruent to
$\Pi_{\sigma \in {\rm Gal}(K/{\bf Q})}
\sigma(\alpha)^{m_{i,\beta,\sigma}}$ mod $p'$ by choice of $\beta$, 
are congruent mod $\wp'$ under some ordering, for $\wp'$ any prime above $p$
in a sufficiently large number field. 
From this and the fact that $p'$ was chosen so that, whenever
$\Pi_\sigma \sigma(\alpha)^{m_{\sigma}}-1$ is not coprime to $p'$
and $|m_{\sigma}| \leq N$, then all the $m_{\sigma}$'s are 0,
the claim follows.
Such elements $\beta$ surject onto $({\cal O}_K/m_{\wp}p{\cal O}_K)^*$ 
for almost all primes $p$ of ${\bf Z}$. This follows from the Cebotarev density theorem as the only condition on the $\beta$'s is that they generate split prime ideals not in $S$ and that for a fixed prime $p'$ of ${\bf Z}$ they be congruent 
to a fixed number $\alpha$ mod $p'$. We denote the common value of 
$m_{i,\beta,\sigma}$ for all such $\beta$'s by $m_{i,\sigma}$. (At this point we know that
the multiplicities with which distinct tuples occur in $\langle (m_{i,r,\sigma})_{\sigma} \rangle$
is independent of $r$.)

From this we conclude (see Section 4.2 of [K] for more details) that
the homomorphisms
$\rho_{\wp}^m:({\cal O}_K/{\sf m}_{\wp}p{\cal O}_K)^* \rightarrow GL_n({\bf F}_{\wp})$ 
for almost all primes $\wp$ are the direct sums 
of the homomorphisms that arise from reducing the homomorphisms
$x \rightarrow \Pi_{\sigma \in {\rm Gal}(K/{\bf Q})}\sigma(x)^{m_{i,\sigma}}$
mod $\wp$. In particular $\rho_{\wp}^m:({\cal O}_K/{\sf m}_{\wp}p{\cal O}_K)^* \rightarrow GL_n({\bf F}_{\wp})$
factors through
$({\cal O}_K/p{\cal O}_K)^*$. 
Thus, remembering that the ${\sf m}_{\wp}$'s were divisible
by only finitely many primes, 
we conclude that $\rho_{\wp}$ factors 
through $({\cal O}_K/{\sf m}p{\cal O}_K)^*$
for a non-zero ideal $\sf m$ independent of $p$.
Then again as in proof of Theorem 1 of [K] (Section 4.2 of [K]) we repeat the argument above using principal split prime ideals that are congruent
to 1 mod $\sf m$. Namely we see that for such prime ideals $r$, the roots of the
characteristic polynomial $f_{r}(X)$ are $\Pi_{\sigma}\sigma(r)^{m''_{i,r,\sigma}}$ with  integers
$m''_{i,r,\sigma}$.
Using the argument above we conclude that the $m''_{i,r,\sigma}$'s 
are independent of $r$, and
thus we can set $m''_{i,r,\sigma}=m''_{i,\sigma}$. This is enough to prove
that the
system of homomorphisms 
$\rho_{\wp}:({\cal O}_K/{\sf m}_{\wp}p{\cal O}_K)^* \rightarrow GL_n({\bf F}_{\wp})$
{\it arise} from algebraic characters
of $K^*$, where $K^*$ is considered as the ${\bf Q}$-valued points 
of the algebraic group ${\bf Res}_{K/{\bf Q}}({\bf G}_m)$ over ${\bf Q}$, 
that have a subgroup of finite index of the units ${\cal O}_K^*$ in their 
kernel. From this we conclude that the compatible system
$(\rho_{\wp})$ of Theorem \ref{main} is
the direct sum of the one-dimensional compatible systems that arise from  Hecke
characters $\chi_i$ for $i=1,\cdots,n$ 
such that $\chi_i$ has infinity type $(m''_{i,\sigma})_{\sigma}$.
This end the proof of the theorem.

\vspace{3mm}

\noindent{\bf Remark:} There is a simpler {\it d{\'e}vissage} argument, in which one takes determinants
of the compatible system and deduces 
Theorem \ref{main} from Theorem 1 of [K], to prove Theorem \ref{main} 
in the case when we know that the compatible 
system we are dealing with is {\it integral} (see Definition 1 of [K]).

\section{Analogs for function fields}

While in both [K] and the present paper we have been 
concerned with representations of absolute Galois groups of number fields,
it looks plausible that the methods of [K] and this paper should carry over to this setting and lead to the 
classification of abelian compatible mod $\wp$ systems of 
absolute Galois groups of function fields. 

Let ${\bf F}$ be a finite field of characteristic $p$ and consider the ring of polynomials
${\bf F}[t]$ and rational functions ${\bf F}(t)$: these will serve as our analogs of ${\bf Z}$ and ${\bf Q}$ respectively. 
Any function field, say $K$, we consider will
be a finite separable extension of ${\bf F}(t)$, 
and the ring of integers ${\cal O}_K$ will
be the integral closure of ${\bf F}[t]$ in $K$. Thus by our choices if $K$ corresponds to the function field
of a smooth projective curve $X$ that is geometrically irreducible over a finite extension of ${\bf F}$,
then $X$ comes equipped with a finite set of infinite places. We denote by $G_K$ the Galois
group of the separable closure $K^s$ of $K$ over $K$.

To get the right analog of compatible systems
in this setting one would need to consider
the field of rationality $L$, which is part 
of definition of compatible systems, 
to be a finite extension of the function field $K$.
In [G] there is an account of compatible systems of $\wp$-adic representations in this 
setting. Here is the definition of compatible mod $\wp$ systems in this setting.
\begin{definition}\label{compatibleFF}
 Let $K$ and $L$ be function fields as above, with a choice of infinite places of $K,L$ as above, and $S,T$ finite sets of places 
 of $K$ and $L$ respectively that contain the infinite places. An $L$-rational 
 strictly compatible system $\{\rho_{\wp}\}$ 
 of $n$-dimensional mod $\wp$ representations of
 $G_K$ with defect set $T$ and ramification set $S$,
 consists of giving for each finite place $\wp$ of $L$ not in $T$ a 
 continuous, semisimple representation $$\rho_{\wp}:G_K \rightarrow GL_n({\bf
 F}_{\wp}),$$ for ${\bf F}_{\wp}$ the residue field of ${\cal O}_L$ at
 $\wp$ of characteristic $p$, that is \begin{itemize} \item  unramified at the places
 outside $S \cup$ $\{$all places of $K$
 above the place of ${\bf F}_q[T]$ below $\wp$ $\}$ \item
 $\rho_{\wp}(G_{\infty_i})$ is trivial where $G_{\infty_i}$ is a decomposition group at 
 the infinite places $\infty_i$ 
 \item for each place 
 $r$ of $K$ not in $S$ there is a monic polynomial $f_r(X) \in L[X]$ such that for all places $\wp$ of
 $L$ not in $T$, and that do not lie above the prime of ${\bf F}[T]$ below $r$, and
 such that $f_r(X)$ has coefficients that are integral at $\wp$, 
 the characteristic polynomial of $\rho_{\wp}({\rm Frob}_r)$
 is the reduction of $f_r(X)$ mod $\wp$, where ${\rm Frob}_r$ is the
 conjugacy class of the Frobenius at $r$ in the Galois
 group of the extension of $K$ that is the fixed field of the kernel of
  $\rho_{\wp}$. \end{itemize} 
\end{definition} 
The condition of being split, or at least potentially split via a base change independent of $\wp$, 
at the infinite places is natural in this context as pointed out to us
by Gebhard B{\"o}ckle. There is a definition of Hecke characters for function fields
in [Gr], that depends on the choice of infinite places, and there it is indicated how to attach compatible system
of $\wp$-adic representations to Hecke characters in this situation that are (potentially) split at the
infinite places. Then we expect that the methods here should be able to prove that
 an $L$-rational strictly compatible system $\{\rho_{\wp}\}$ 
 of abelian, semisimple mod $\wp$ representations of
 $G_K$ is a direct sum of one-dimensional 
 compatible systems
 each of which arises from a Hecke character. From this will follow the following 2 statements:
 
1. An $L$-rational strictly compatible system $\{\rho_{\wp}\}$ 
 of abelian, semisimple mod $\wp$ representations of
 $G_K$ lifts to a compatible system of $n$-dimensional 
  $\wp$-adic representations
  $\{\rho_{\wp,\infty}\}$ (split at all infinite places).
 
2. An $L$-rational strictly compatible system of abelian, semisimple
$n$-dimensional 
  $\wp$-adic representations
  $\{\rho_{\wp,\infty}\}$ of $G_K$
  arises from Hecke characters.

It will be useful to have a systematic exposition that will give details of
construction of Galois representations associated to Hecke characters in this setting (generalising and giving
more details the work in [Gr]) 
and also prove the statements above. This will lead to a proof that the compatible system of 
characteristic $p$ Galois representations attached by G.~B{\"o}ckle,
R.~Pink and others to Drinfeld modular forms arise from Hecke characters.
It is also seems plausible that the second statement above can be strengthened 
considerably: 

\vspace{3mm}

\noindent{\bf Question:} If we have a (continuous) abelian semisimple representation $\rho: G_K \rightarrow GL_n(L_{\wp})$
that is (potentially) split at the infinite places, with $\wp$ a finite place of a function field $L$ and $L_{\wp}$ the completion of $L$ at $\wp$, 
unramified oustide a finite set of places $S$ that includes the infinite
places, and for places $v$ not in $S$, the characteristic polynomial of $\rho({\rm Frob}_v)$ has coefficients in $L$, then 
is $\rho$ the sum of 1-dimensional compatible systems that arise
from Hecke characters in the sense of Gross in [Gr]? 

\vspace{3mm}

\noindent The analog of this question for number fields is answered affirmatively by G.~Henniart in [He]  
using results of Waldschmidt in transcendental
number theory. Of course a question like the one above cannot be answered using only
the elementary algebraic techniques of the present paper.

\vspace{3mm}

\noindent{\bf Acknowledgements:} The author would like
to acknowledge helpful conversations and correspondence with
Srinath Baba, Gebhard B{\"o}ckle, Guy Henniart, Dinesh Thakur
and Michel Waldschmidt.

\section{References}

\noindent [CS] Corrales, C., and Schoof, R., {\it The support
problem and its elliptic analogue}, J. of Number Theory 64 (1997), 276--290.

\vspace{3mm}

\noindent [G] Goss, D., {\it Basic structures of function field arithmetic}, Springer-Verlag, 1996.

\vspace{3mm}

\noindent [Gr] Gross, B., {\it Algebraic Hecke characters for function fields},
in S{\'e}minaire  de Th{\'e}orie des
Nombres, Progress in Math. 22 (1982), 87--90, Birkhauser.

\vspace{3mm} 

\noindent [He] Henniart, G., {\it Repr{\'e}sentations $\ell$-adiques 
ab{\'e}liennes}, in S{\'e}minaire  de Th{\'e}orie des
Nombres, Progress in Math. 22 (1982), 107--126, Birkhauser.

\vspace{3mm}

\noindent [K] Khare, C., {\it Compatible systems of mod $p$ Galois representations and Hecke characters},
Mathematical Research Letters 10 (2003), vol. 1, 71--84.

\vspace{3mm}

\noindent [S] Serre, J-P., {\it Abelian $\ell$-adic representations
and elliptic curves}, Addison-Wesley, 1989.





\vspace{3mm}

\noindent Dept of Math, University of Utah,
155 S 1400 E, Salt Lake City, UT 84112, USA. e-mail: shekhar@math.utah.edu

\noindent School of Mathematics, TIFR, Homi Bhabha Road, Mumbai 400 005,
INDIA. e-mail: shekhar@math.tifr.res.in

\end{document}